\documentclass[a4paper]{article}
\usepackage{typearea}
\areaset[1cm]{158mm}{234mm}
\let\production 0

\usepackage{xcolor}
\usepackage{graphicx}
\usepackage{psfrag}
\usepackage{caption}
\usepackage{subcaption}

\usepackage{fancybox}
\usepackage{fancyvrb}

\usepackage{amsmath}
\usepackage{amsthm}
\usepackage{amssymb}

\usepackage{latexsym,amssymb}
\usepackage{url}
\usepackage{amsfonts}
\usepackage{hhtensor}

\usepackage{indentfirst}

\def\RR{\mathbb{R}}

\def\DD{\mathbb{D}}
\def\HH{\mathbb{H}}
\def\LL{\mathbb{L}}

\def\NN{\mathbb{N}}
\def\TT{\mathbb{T}}

\def\G{{\cal G}}

\def\M{{\cal M}}

\def\F{{\cal F}}

\def\dfrac{\displaystyle\frac}
\def\tens{\mathbf}
\def\tens{\boldsymbol}
\def\media_#1{\mathchoice{\mathop{\hskip2pt{-}\hskip-10.3pt\int\hskip-4pt}%
    \nolimits_{#1}}%
    {\mathop{-\hskip -8.8pt\int}\nolimits_{#1}}{}{}}

\let\epsi\varepsilon
\let\ffi\varphi

\def\div{\mathrm{div}}
\def\Ker{\mathrm{Ker}}
\def\Im{\mathrm{Im}}
\def\dint{\displaystyle\int}

\theoremstyle{plain}
\newtheorem{theorem}{Theorem}

\newtheorem{lemma}[theorem]{Lemma}

\theoremstyle{definition}

\newtheorem{remark}[theorem]{Remark}

\title{Donati representation theorem for periodic functions in relation to homogenization theory}
\author{Cristian Barbarosie$^1$\thanks{\texttt{cabarbarosie@fc.ul.pt}} \and Anca-Maria Toader$^2$\thanks{\texttt{atoader@fc.ul.pt}, 
both authors are from DM-FCUL and CMAFcIO Centro de Matem\'atica, Aplica\c c\~oes Fundamentais e Investiga\c c\~ao Operacional, 
Departamento de Matem\'atica da Faculdade de Ci\^encias, 
Universidade de Lisboa, 1749-016 Lisboa, Portugal.
}}
\date{\today}

\begin{document}
	\maketitle
	
	\begin{abstract}
This paper discusses properties of periodic functions, focusing on (systems of) partial differential
equations with periodicity boundary conditions, called {\sl cellular problems}.
These cellular problems arise naturally from the asymptotic study of PDEs with rapidly oscillating
coefficients, called {\sl homogenization theory}.
We believe the present paper may shed a new light on well-known concepts, for instance by showing
hidden links between Green's formula, the div-curl lemma, compensated compactness and Donati's 
representation theorem.
We state and prove three extensions of Donati’s Theorem adapted to the periodic framework which,
beyond their own importance, are essential for understanding the variational formulations of
cellular problems in strain, in stress and in displacement.
Section \ref{sec:traces} presents a self-contained study of properties of
traces of a function and their relations with periodicity properties of that function.

	\end{abstract}



  


\section{Introduction and historical considerations}
\label{sec:history}

This paper presents a study of properties of periodic functions, vector fields and matrix fields,
in relation with (systems of) partial differential equations arising from homogenization theory.
It is closely related to another paper by the same authors \cite{BT2022-B} which can be
viewed as a continuation of this one.
We believe these two papers may shed a new light on well-known concepts like Green's formula,
traces of functions, variational formulation of cellular problems and Donati's representation theorem,
in the context of periodic homogenization.

Section \ref{sec:history} presents some historical facts about representation theorems.
Section \ref{sec:homogenization} recalls concepts and results from homogenization theory.
Section \ref{sec:Korn-Green} recalls two fundamental results, adapted to periodic functions.
Although its main goal is to prove Lemmas \ref{lemma:traces-5} and \ref{lemma:traces-7},
section \ref{sec:traces} can be viewed as a self-contained study of properties of traces
in the periodic context.
In section \ref{sec:compens-compact} we highlight links between Green's theorem and Tartar's
div-curl Lemma.
Section \ref{sec:LP-T-Donati} is devoted to the main results of this paper\,:
Theorems \ref{thrm:Donati.1}, \ref{thrm:Donati.2} and \ref{thrm:Donati.3} are
three extentions of Donati's theorem for functions satisfying periodicity boundary conditions.

We start by some historical facts in order to point out the importance of Donati's Theorem and
to place in time the development of its extensions.

It is well known that the gradient of a scalar function has zero curl.
The reverse is also true, at least for a simply-connected domain\,: any irrotational
vector field is the gradient of some function.
In differential geometry, the Poincar\'e lemma generalizes this property by means of differential forms.

One may ask a similar question about matrix fields\,: when is a matrix field the gradient of a
vector-valued function ?
The answer to this question is straightforward because the problem can be split\,:
each row of the matrix can be viewed as the gradient of a scalar function.
But if we slightly change the question by asking when is a matrix field the {\em symmetric part} of
the gradient of some vector-valued function, things get more complicated.
This question cannot be answered by analysing separately each row of the matrix.

The symmetric part of the gradient of a vector function arises naturally in the context of
continuum mechanics\,: in infinitesimal strain theory, the strain field is the symmetric part
of the displacement gradient.

In 1864 Saint Venant announced a theorem which would later bear his name and
which was proved rigorously by E.\ Beltrami in 1886\,:

\begin{theorem} [Saint Venant]  \label{thrm:Saint_Venant}
 Assume that the open set \,$\Omega$ is simply-connected. Then there exists a vector field
$v\in C^3(\Omega)$ such that $\tens e = \nabla_s v$ in $\Omega$ if (and clearly only if, even if $\Omega$ is not simply-connected) 
the functions $\tens e_{ij}$ are in
the space $C^2(\Omega)$ and they satisfy
$$
R_{ijkl}(\tens e) := \partial_{lj} \tens e_{ik} + \partial_{ki} \tens e_{jl} -\partial_{li} \tens e_{jk} -\partial_{kj} \tens e_{il} = 
0 \hbox{ in }\Omega \hbox{ for all } i, j, k,l \in \{1, 2, 3\}.
$$
\end{theorem}

As we can see, the compatibility conditions for a matrix field to be a symmetric gradient
involve second-order derivatives,
which is very different from the ``curl-free'' condition for gradients.
 
In 1890 Donati gave a classic result\,: a necessary condition on a strain tensor 
field to be orthogonal to all divergence free stress tensor fields. The
condition is written in terms of second order derivatives of the strain tensor
field, and for simply connected domains is equivalent to the so-called Saint
Venant compatibility condition.

\begin{theorem} [Donati]
\label{them:Donati}
If $\Omega$ is an open subset of $\RR^3$ and the components $\tens e_{ij}$ of a symmetric matrix field
$\tens e = (\tens e_{ij} )$ are in the space $C^2(\Omega)$ and they satisfy:
$$
\int_\Omega \tens e_{ij} \tens s_{ij} dx = 0 \hbox{ for all } \tens s = (\tens s_{ij} ) 
\in \DD_s(\Omega) \hbox{ such that } \div \tens s = 0 \hbox{ in } \Omega,
$$
where $\DD_s(\Omega)$ denotes the space of all symmetric tensor fields whose components are infinitely differentiable in $\Omega$
and have compact supports in $\Omega$, then
$$
CURL \, CURL \,\tens e = 0 \hbox{ in } \Omega,
$$
where $ (CURL \tens e)_{ij} = \epsi_{ilk} \partial_l \tens e_{jk}$ for any $\tens e \in \DD^\prime (\Omega)$.
\end{theorem}

Donati’s Theorem thus provides, once combined with Saint Venant’s Theorem,
another characterization of symmetric matrix fields as linearized strain tensor fields,
at least for simply-connected open subsets $\Omega$ of $\RR^3$.

An extension of Donati’s Theorem was obtained by T.W.\ Ting 84 years later, in 1974.
In \cite{Ting} Ting stated a necessary condition on a strain tensor field to be orthogonal
to all divergence free stress tensor fields, in terms of Lebesgue and Sobolev spaces.

\begin{theorem} [Ting]
\label{thrm:Ting}
{\it If $\Omega$ is bounded and has a Lipschitz continuous boundary and if the components $\tens e_{ij}$ of the symmetric tensor field 
$\tens e$ are in $L^2(\Omega)$ and again satisfy:
$$
\dint_\Omega \tens e_{ij} \tens s_{ij} dx = 0 \hbox{ for all } \tens s = (s_{ij} ) \in  \DD_s(\Omega)
\hbox{ such that } \div \tens s = 0 \hbox{ in }\Omega ,
$$
then there exists $v \in (H^1(\Omega))^3$ such that $\tens e = \nabla_s v$ in $\LL^2_s(\Omega)$.}
\end{theorem}

In 2006, C.\ Amrouche, P.G.\ Ciarlet, L.\ Gratie and S.\ Kesavan, see \cite{Amrouche2006},
proved that the condition in Donati's Theorem is equivalent to Saint Venant's compatibility conditions.

\begin{theorem} [Amrouche, Ciarlet, Gratie \& Kesavan]
\label{thrm:SV_compat_cond_equiv_curl_curl_eq_zero}
The Saint Venant compatibility conditions $R_{ijkl}(\tens e) = 0$ in $\Omega$ are equivalent to
the relations
$CURL\, CURL\, \tens e = 0$ in $\Omega$.
\end{theorem}

In the last 20 years P.\ G.\ Ciarlet, P.\ Ciarlet, C.\ Amrouche and their collaborators wrote several papers
where extensions of the Donati’s theorem were
obtained and applied to PDEs with Dirichlet, Neumann and mixed boundary conditions, see, for example 
\cite{Ciarlet-Ciarlet2005}, \cite{Amrouche2006}, \cite{Ciarlet-G-K2011}.

However, the periodicity condition was not considered and remained open until now.
In Section \ref{sec:LP-T-Donati} we state and prove a series of extensions on Donati's Theorem
adapted to the periodic context, in relation to the theory of homogenization.

\section{The homogenization theory as motivation}
\label{sec:homogenization}

We begin by describing the general periodicity notion that we have been using since \cite{BT2010}
in the homogenization theory.

Consider a lattice of vectors in $ \RR^3 $ (an additive group $ \G $ generated by three
linearly independent vectors $ \vec g_1, \vec g_2, \vec g_3 $).
Define the parallelepiped  $ Y \subset \mathbb{R}^3 $\,:
$$
Y = \{ s_1 \vec g_1 + s_1 \vec g_2 + s_1 \vec g_3 : s_1, s_2, s_3 \in [0, 1]\}.
$$
$Y$ is called ``periodicity cell''.

A function $\ffi $ is said to be $\G$-periodic (or $Y$-periodic) if it is invariant to translations
with vectors in $ \G $.
In a periodic microstructure the rigidity is a fourth order tensor field
$ \tens C : \mathbb{R}^3 \mapsto \mathbb{R}^{81} $ which is $\G$-periodic, 
that is, it varies according to a periodic pattern.
Denote by $ \tens B $ the compliance tensor, that is, the inverse tensor of $ \tens C $.

According to the homogenization theory, the macroscopic behaviour of a body with periodic structure
is described by a constant homogenized elastic tensor $ \tens C^H $.
It is possible to define $ \tens C^H $ in terms of the solutions of the cellular problems
which are PDEs subject to periodicity conditions:


\begin{equation}  \label{eq:cell-pb}
\left\{\begin{matrix}
\div (\tens C\, \tens e(u_A))  =  0 \hbox{ in } \RR^3 \\
u_A(y)  =  A y + \ffi_{ A}(y), 
\quad \hbox{ with }  \ffi_A \quad \G-\hbox{periodic} ,
\end{matrix}\right.
\end{equation}
where $A$ is a given macroscopic strain (a $ 3 \times 3 $ symmetric matrix) and $\tens e$
represents the corresponding microscopic strain, that is, the symmetric part of the gradient
of $ u_A $.
In the sequel, for historical reasons, we shall use also the notation $ \nabla\!{}_s $ for the
symmetric part of the gradient, when the intention is to focus on the displacement $ u_A $.

The solution $ u_A $ depends linearly on the matrix $A$ and verifies 
\begin{equation}
\label{eq:A}
A = \media_Y  \tens e  (u_A) = \frac{1}{|Y|} \int_Y  \tens e (u_A) \,.
\end{equation}

\noindent The macroscopic stress associated to $ u_A $ is defined by
\begin{equation}  \label{eq:S}
S = \media_Y \tens C\, \tens e  (u_A) = \frac{1}{|Y|} \int_Y \tens C\, \tens e (u_A) 
\end{equation}
and consequently depends linearly on $A$.
The symbol $ \displaystyle\media_Y $ denotes the average value of the integrand
on the periodicity cell $Y$.
	
The homogenized elastic tensor is then defined, for all macroscopic strains $A$,
through the equality $\tens C^H A = S $, that is,
\begin{equation} \label{CHA}
\tens C^H A = \media_Y \tens C\, \tens e  (u_A) .
\end{equation}

\noindent $\tens C^H$ can be also defined, equivalently, in terms of energy products\,:
\begin{equation} \label{CHAB}
\langle \tens C^H A, B \rangle =
\media_Y \langle \tens C\, \tens e (u_A), \tens e (u_B)\rangle,
\end{equation}
where $A$ and $B$ are any two symmetric matrices (strains).
See Remark \ref{rem:hom-tens-energ-prod} below.

We present in the sequel some notations used throughout the paper.

Denote by $ L^2_\# (\RR^3,\RR^3) $ the space of vector fields in $ L^2_{\mbox{\small loc}} (\RR^3,\RR^3) $
which are $ \G $-periodic, endowed with the $ L^2 (Y,\RR^3) $ norm.
Denote by $ H^1_\# (\RR^3,\RR^3) $ the space of vector fields in $ H^1_{\mbox{\small loc}} (\RR^3,\RR^3) $
which are $ \G $-periodic, endowed with the $ H^1 (Y,\RR^3) $ norm.
Equivalently, $ H^1_\# (\RR^3,\RR^3) $ can be viewed as the completion
in the norm of $ H^1(Y,\RR^3) $ of the space of functions in $ C^\infty(\RR^3, \RR^3) $ which are 
$\G$-periodic. Denote by $ H^{-1}_\# (\RR^3, \RR^3) $ the dual space of $ H^1_\# (\RR^3, \RR^3) $.

For the sake of abreviation we shall use the notations $ L^2_\# $ for $ L^2_\# (\RR^3,\RR^3) $,
$H^1_\# $ for $H^1_\# (\RR^3,\RR^3)$ and $ H^{-1}_\# $ for $ H^{-1}_\# (\RR^3, \RR^3) $.

Denote by $ \LL^2_\# $ the space of matrix fields in $ L^2_{\mbox{\small loc}} (\RR^3,\RR^9) $
which are $ \G $-periodic, endowed with the $ L^2 (Y,\RR^9) $ norm;
we add the subscript $s$ for symmetric matrices\,: $ \LL^2_{\# s} $.
Denote by $ \HH^1_\# $ the space of matrix fields in $ H^1_{\mbox{\small loc}} (\RR^3,\RR^9) $
which are $\G$-periodic, endowed with the $ H^1 (Y,\RR^9) $ norm;
we add the subscript $s$ for symmetric matrices\,: $ \HH^1_{\# s} $.

Denote by $\HH_{\# s} (\div )$ the space

$$
\HH_{\# s} (\div ) = \bigl\{ \tens\mu \in \LL^2_{\# s} \mid \div \tens\mu \in  L^2_\# \bigr\}\,,
$$

\noindent where, with Einstein's repeated index notation, the operator
$ \div : L^2_{\mbox{\small loc}} (\RR^3,\RR^9) \to H^{-1}(\RR^3, \RR^3) $ is defined by
$ \div \tens\mu = (\mu_{ij,j})_{1\le i\le 3} $.

Consider $ LP $ the space of linear plus periodic displacements defined in $\RR^3$\,:
$$
\begin{matrix}
LP = \bigl\{ u : \RR^3 \to \RR^3 \mid u(y)  = Ay + \ffi (y), \\
 A \in \M_3^s(\RR), \ffi \in H^1_\# \bigr\},
\end{matrix} 
$$
where $ \M_3^s(\RR) $ denotes the space of symmetric $ 3\times 3 $ matrices with real elements.
$ LP $ is a Hilbert
space since it is a direct sum between a finite dimensional space (of all linear maps) and $ H^1_\# $.
$ LP $ can equivalently be viewed as the space of all vector fields $ u \in H^1_{\mbox{\small loc}}
( \RR^3 ) $ whose derivatives are $\G$-periodic.

For an arbitrarily fixed strain matrix $ A\in  \M_3^s(\RR) $, denote by $ LP(A) $ the set of linear plus periodic 
displacements having the linear part $ Ay $ \,:
$$
\begin{matrix}
  LP(A)= \bigl\{ u : \RR^3 \to \RR^3 \mid u(y) = Ay + \ffi (y), \ffi \in H^1_\# \bigr\}
  = \{ Ay \} + H^1_\# .
\end{matrix} $$
Thus the last equation in (\ref{eq:cell-pb}) is equivalent to $ u_A\in LP(A) $.
Note that $ LP(0) = H^1_\#  $ is a closed subspace of $ LP $ and for
a given strain matrix $ A\in \M_3^s(\RR) $ the set $ LP(A) $ is a translation of $ LP(0) $. 
Moreover,
$$ LP = \displaystyle\bigcup_{A\in \M_3^s(\RR)} LP(A). $$ 

Define the space of linear plus periodic functions having zero average by\,: 
\begin{equation}
\label{eq:LP_0}
LP_0 :=  \{ v \in LP \mid \media_Y v dy =0 \}.
\end{equation}

For a given $ 3\times 3 $ symmetric matrix $ S \in \M_3^s(\RR) $, define the space of stresses 
\begin{equation} \label{eq:def-TT-S}
\TT(S) = \{ \tens \mu \in \HH_{\# s}(\div) \mid \div \tens \mu = 0, \quad \media_Y\tens \mu\, dy = S  \}.
\end{equation}

Consider the space $ \TT $ defined as
\begin{equation} \label{eq:def-TT}
\TT : = \bigcup_{S\in \M_3^s(\RR)} \TT(S) =
\{ \tens \mu \in \HH_{\# s} (\div ) \mid \div \tens \mu = 0 \hbox{ in } L^2_\# \} . 
\end{equation}

\noindent By $ \TT(0) $ we denote the space of stresses having zero mean.
Define the orthogonal complement of $ \TT(0) $\,:
\begin{equation} \label{eq:def-TT-0-ort}
\TT(0)^\perp  : = \bigl\{ \tens e \in \LL^2_{\# s}  \mid \dint_Y \langle \tens e, \mu \rangle =0,
\forall \mu \in  \TT(0) \bigr\} . 
\end{equation}

The variational formulation of the cellular problem (\ref{eq:cell-pb}), when the macroscopic strain
$ A \in \M^s_3(\RR) $ is given, has the following form\,:

\begin{equation} \label{eq:var-form-displ-macrostrain}
\left\{
\begin{matrix} \mbox{find } u_A \in LP_0 \cap LP(A) \mbox{ such that } \\
  \displaystyle\int_{Y} \langle \tens C\, \tens e(u_A) ,\tens e(v)\rangle\, dy =0
  \ \ \forall v \in H^1_\# . 
\end{matrix} \right. 
\end{equation}

Several other variational formulations of the same cellular problem (\ref{eq:cell-pb}) are stated
in Section 2 of \cite{BT2022-B}.

\section{Korn inequality, Green's formula}
\label{sec:Korn-Green}

The current Section deals with two fundamental results adapted to the periodicity framework.

\begin{theorem}
\label{thrm:Korn}
The Korn inequality below holds for a positive constant $C$ and for all $v$ in $H^1_{\# }$
\begin{equation}
\label{eq:Korn_ineq}
\|v\|_{H^1_\#} \le C( \| v \|_{\LL^2_\# } + \| \nabla\!{}_s v \|_{\LL^2_\# })
\end{equation}
\end{theorem}

\begin{proof}
  Since any function of $H^1_{\# }$ has its restriction to $Y$ belonging to $H(Y)$, the Korn inequality
  above is a direct consequence of the Theorem 2.1 in \cite{Ciarlet-Ciarlet2005}.
  However, to turn this work more self contained we give a proof which is analogous to the one given in
  \cite{Ciarlet-Ciarlet2005} for the Theorem 2.1 and goes back to Theorem 3.2, Chap.\ 3 of
  \cite{Duvaut-Lions}.
We prove that the space  $H^1_{\# }$ coincides with 
$$
K_\# := \{ v \in L^2_{\#} : \nabla\!{}_s v \in \LL^2_{\# s} \}.
$$
The inclusion $H^1_{\# } \subset K_\# $ is obvious. To prove the other inclusion, consider $v \in K_\#$. Then $v\in L^2_{\#}=L^2(Y)$ and
consequently the derivatives $\partial_k v_i \in H^{-1}(Y)$ and hence
$\partial_{jk}v_i = \partial_j e_{ik}(v) + \partial_k e_{ij}(v) - \partial_i e_{jk}(v)$ belongs to $H^{-1}(Y)$.
Then a fundamental Lemma by J.\ L.\ Lions ($v\in H^{-1}(\Omega), \partial v \in H^{-1}(\Omega) \Rightarrow v \in L^2(\Omega)$ ) implies that
$\partial_k v_i$ belongs to $\LL^2_\#$, which completes the proof of the equality between the two spaces. 
The Korn inequality results as a consequence of the closed graph theorem applied to the identity mapping from $H^1_{\# }$ into $K_\#$ which 
turns out to be surjective and continuous. 

\end{proof}

\begin{theorem} \label{thrm:Green}
The Green's formula 
\begin{equation}  \label{eq:Green}
\dint_Y \langle \tens \mu , \nabla\!{}_s v \rangle\, dy + \dint_Y (\div \tens \mu) \cdot v\, dy = 0
\end{equation}
holds for all $\tens \mu \in \HH_{\# s} (\div) $ and all $v \in H^1_\# (\RR^3,\RR^3)$.		
\end{theorem}

The proof of Green's formula above is based on the notion of trace and its properties in the context
of periodic functions, presented in section \ref{sec:traces}.

\begin{proof}
  Since $ \tens\mu $ is symmetric, we can replace, in the first integral, $ \nabla\!{}_s v $ by $ \nabla v $.
  Corollary 2.1 in \cite{GiraultRaviart} states that
  $$
  \dint_Y \tens\mu_{ij} v_{i,j}\, dy + \dint_Y  \tens\mu_{ij,j} v_i\, dy =
  \int_{\partial Y} \tens\mu_{ij} n_j v_i\,d\sigma
  $$
  Th\'eor\`eme 1.1 in \cite[Chapitre 3]{Necas1967} also states the above equality,
  under somewhat stronger assumptions.
  We write $ \partial Y $ as the union of 6 faces (we are in three dimensions)
  $ \partial Y = F_1^+ \cup F_1^- \cup F_2^+ \cup F_2^- \cup F_3^+ \cup F_3^- $.
  Thus
  $$
  \int_{\partial Y} \tens\mu_{ij} n_j v_i\,d\sigma =
  \sum_k \Bigl( \int_{F_k^+} \mu_{ij} n_j v_i\,d\sigma +
  \int_{F_k^-} \tens\mu_{ij} n_j v_i\,d\sigma \Bigr)
  $$
  Applying Lemmas \ref{lemma:traces-5} and \ref{lemma:traces-7} together with Remark \ref{rem:traces-1}
  in section \ref{sec:traces} below, we conclude that the above sum is zero.
\end{proof}

\section{Traces of functions in $ H^1 $}
\label{sec:traces}

This section presents a self-contained study of traces of functions and their properties.
The first half of this section dwells on the general setting with no periodicity assumptions.
The second half is about periodic functions and their traces on the boundary of the periodicity cell.
Lemmas \ref{lemma:traces-5} and \ref{lemma:traces-7} are useful for proving Green's formula
(Theorem \ref{thrm:Green} in Section \ref{sec:Korn-Green}).

In this section we show drawings and use notations relative to two-dimensional domains
(functions of two variables).
However, it is easy to see that the results here presented are valid in any space dimension.
Also, we refer to scalar-valued functions $f\in H^1$ and vector fields $ \mu\in H(\div) $ but it is clear
that similar results hold for vector fields $u\in H^1$ and matrix fields $ \tens\mu\in \HH(\div) $.

It is well-known that a function in $ H^1(\Omega) $ has a trace on $ \partial\Omega $.
Since any function in $ H^1_\# $ can be restricted to a function in $ H^1(Y) $,
it will have a trace on $ \partial Y $.
However, new questions arise because a function in $ H^1_\# $ is defined on both sides
of $ \partial Y $.
Thus, it has one trace from the interior of $Y$ and one trace from ``outside'', which can be viewed
as coming from a neighbour translation of $Y$.
The present section dwells on these details.

\begin{figure}[ht] \centering
  \psfrag{O1}{$ \Omega_1 $}
  \psfrag{O2}{$ \Omega_2 $}
  \psfrag{G}{$ \Gamma $}
  \psfrag{n1}{$ n^{(1)} $}
  \psfrag{n2}{$ n^{(2)} $}
  \includegraphics[width=55mm]{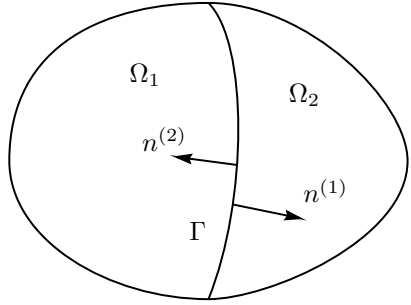}
  \caption{Two domains with a common interface}
  \label{fig duas batatas}
\end{figure}

First, let us forget momentarily about the periodic setting;
consider two domains $ \Omega_1 $ and $ \Omega_2 $ having a common piece of boundary $ \Gamma $,
as shown in Figure \ref{fig duas batatas}.
Denote by $ \Omega $ their union $ \Omega_1 \cup \Omega_2 \cup \Gamma $.

\begin{lemma}\label{lemma:traces-1}
  Consider a function $ u \in H^1 ( \Omega ) $.
  Then the trace $ tr^{(1)}_\Gamma u $ of $u$ on $ \Gamma $ coming from $ \Omega_1 $
  (that is, the trace of $ u\!\!\mid_{\Omega_1} $ on $ \Gamma $) is equal to the trace
  $ tr^{(2)}_\Gamma u $ of $u$ on $ \Gamma $
  coming from $ \Omega_2 $ (that is, the trace of $ u\!\!\mid_{\Omega_2} $ on $ \Gamma $).
\end{lemma}

\begin{proof}
  Let $ \varphi $ be a $ C^\infty $ function whose support is compact in $ \Omega $
  (thus, it vanishes on $ \partial\Omega $ but may be non-zero on $ \Gamma $).
	Green's formula, applied separately in $ \Omega_1 $ then in $ \Omega_2 $, states that
  \begin{align*}
    \int_{\Omega_1} \partial_i u\, \varphi \,dx =
    \int_\Gamma tr^{(1)}_\Gamma u\, \varphi\, n^{(1)}_i\, d\sigma -
    \int_{\Omega_1} u\, \partial_i \varphi\, dx \\
    \int_{\Omega_2} \partial_i u\, \varphi \,dx =
    \int_\Gamma tr^{(2)}_\Gamma u\, \varphi\, n^{(2)}_i\, d\sigma -
    \int_{\Omega_2} u\, \partial_i \varphi\, dx 
  \end{align*}
  where $ n^{(1)} $ it the unit vector normal to $ \Gamma $ pointing from $ \Omega_1 $ towards
  $ \Omega_2 $ and $ n^{(2)} = - n^{(1)} $.

  On the other hand, we can apply Green's formula in $ \Omega $ as a whole\,:
  \begin{equation*}
    \int_\Omega \partial_i u\, \varphi \,dx = -\int_{\Omega} u\, \partial_i \varphi\, dx
  \end{equation*}

  The above three equalities imply
  \begin{equation*}
    \int_\Gamma ( tr^{(2)}_\Gamma u - tr^{(1)}_\Gamma u )\, \varphi\, n^{(1)}_i\, d\sigma = 0
  \end{equation*}
  for all indices $i$, which means that $ tr^{(2)}_\Gamma u - tr^{(1)}_\Gamma u = 0 $.
\end{proof}

The reverse result also holds\,:

\begin{lemma}\label{lemma:traces-2}
  Consider two functions $ u^{(1)} \in H^1 ( \Omega_1 ) $ and $ u^{(2)} \in H^1 ( \Omega_2 ) $
  such that $ tr^{(1)}_\Gamma u^{(1)} = tr^{(2)}_\Gamma u^{(2)} $.
  Then there exists a function $ u \in H^1 ( \Omega ) $ such that $ u^{(1)} $ is the restriction
  of $u$ to $ \Omega_1 $ and $ u^{(2)} $ is the restriction of $u$ to $ \Omega_2 $.
\end{lemma}

\begin{proof}
  Define $u$ simply by joining the two functions $ u^{(1)} $ and $ u^{(2)} $.
  Since the Lebesgue measure of $ \Gamma $ is zero, the values of $u$ on $ \Gamma $ do not matter.
  It is clear that $ u \in L^2 (\Omega) $; we only need to prove that its first order derivatives
  belong also to $ L^2 (\Omega) $.
  For each index $i$, define $ w_i \in L^2 (\Omega) $ simply by joining the two functions
  $ \partial_i u^{(1)} \in L^2 (\Omega_1) $ and $ \partial_i u^{(2)} \in L^2 (\Omega_2) $.
  It is clear that $ w_i \in L^2 (\Omega) $; we only need to prove that $ \partial_i u = w_i $.
  
  Let $ \varphi $ be a $ C^\infty $ function whose support is compact in $ \Omega $.
  By applying Green's formula separately in $ \Omega_1 $ and in $ \Omega_2 $ like in the proof of
  Lemma \ref{lemma:traces-1} and by noting that the integrals on $ \Gamma $ cancel out, we conclude that
  \begin{equation*}
    \int_\Omega w_i\, \varphi \,dx = -\int_{\Omega} u\, \partial_i \varphi\, dx
  \end{equation*}
  which means that $ \partial_i u = w_i $.
\end{proof}

Similar results can be proven for the space
$$ H (\div, \Omega) = \bigl\{ \mu\in L^2(\Omega) : \div\,\mu \in L^2(\Omega) \bigr\}\,. $$
It is well known that if $ \mu $ is a vector field, $ \mu \in H(\div, \Omega) $,
then the normal component of $ \mu $ has a trace on $ \partial\Omega $.
When we look at a domain which is the union of two smaller domains like in Figure \ref{fig duas batatas},
we see there are two traces on $ \Gamma $, one coming from $ \Omega_1 $ and the other one coming
from $ \Omega_2 $.

\begin{lemma}\label{lemma:traces-3}
  Consider a vector field $ \mu \in H ( \div, \Omega ) $.
  Then the trace $ tr^{(1)}_\Gamma (un^{(1)}) $ of $ \mu n^{(1)} $ coming from $ \Omega_1 $
  (that is, the trace on $ \Gamma $ of the normal component of $ \mu\!\!\mid_{\Omega_1} $)
  is equal to minus the trace $ tr^{(2)}_\Gamma (\mu n^{(2)}) $ of $ \mu n^{(2)} $
  coming from $ \Omega_2 $ (that is, the trace on $ \Gamma $
  of the normal component of $ \mu\!\!\mid_{\Omega_2} $).
  The change in sign is due to the different meanings of the normal $n$, that is, $ n^{(2)} = -n^{(1)} $.
\end{lemma}

\begin{proof}
  Let $ \varphi $ be a $ C^\infty $ function whose support is compact in $ \Omega $
  (this means it vanishes on $ \partial\Omega $ but may be non-zero on $ \Gamma $).
	Green's formula, applied separately in $ \Omega_1 $ then in $ \Omega_2 $, states that
  \begin{align*}
    \int_{\Omega_1} \partial_i \mu_i\, \varphi \,dx =
    \int_\Gamma tr^{(1)}_\Gamma (\mu_i\, n^{(1)}_i)\, \varphi\, d\sigma -
    \int_{\Omega_1} \mu_i\, \partial_i \varphi\, dx \\
    \int_{\Omega_2} \partial_i \mu_i\, \varphi \,dx =
    \int_\Gamma tr^{(2)}_\Gamma (\mu_i\, n^{(2)}_i)\, \varphi\, d\sigma -
    \int_{\Omega_2} \mu_i\, \partial_i \varphi\, dx 
  \end{align*}
  In the above, an implicit sum is assumed in $i$;
  thus, $ \partial_i \mu_i $ means $ \div \mu $ and
  $ \mu_i\, n_i $ is the normal component of $ \mu $.
  
  On the other hand, we can apply Green's formula in $ \Omega $ as a whole\,:
  \begin{equation*}
    \int_\Omega \partial_i \mu_i\, \varphi \,dx = -\int_{\Omega} \mu_i\, \partial_i \varphi\, dx
  \end{equation*}

  The above three equalities imply
  \begin{equation*}
    \int_\Gamma \bigl( tr^{(2)}_\Gamma (\mu_i\, n^{(2)}_i) +
    tr^{(1)}_\Gamma (\mu_i\, n^{(1)}_i) \bigr)\, \varphi\, d\sigma = 0
  \end{equation*}
  which means that $ tr^{(2)}_\Gamma (\mu_i\, n^{(2)}_i) +
  tr^{(1)}_\Gamma (\mu_i\, n^{(2)}_i) = 0 $.
\end{proof}
  
\begin{lemma}\label{lemma:traces-4}
  Consider two vector fields $ \mu^{(1)} \in H ( \div, \Omega_1 ) $ and
  $ \mu^{(2)} \in H ( \div, \Omega_2 ) $ such that $ tr^{(1)}_\Gamma (\mu^{(1)} n^{(1)}) =
  -tr^{(2)}_\Gamma (\mu^{(2)} n^{(2)}) $.
  Then there exists a vector field $ \mu \in H ( \div, \Omega ) $ such that $ \mu^{(1)} $
  is the restriction of $ \mu $ to $ \Omega_1 $ and $ \mu^{(2)} $ is the restriction of
  $ \mu $ to $ \Omega_2 $.
\end{lemma}

\begin{proof}
  Define $ \mu $ simply by joining the two fields $ \mu^{(1)} $ and $ \mu^{(2)} $.
  Since the Lebesgue measure of $ \Gamma $ is zero, the values of $ \mu $
  on $ \Gamma $ do not matter.
  It is clear that $ \mu \in L^2 (\Omega) $; we only need to prove that its divergence
  belongs also to $ L^2 (\Omega) $.
  Define $ w \in L^2 (\Omega) $ simply by joining the two functions
  $ \div \mu^{(1)} \in L^2 (\Omega_1) $ and $ \div \mu^{(2)} \in L^2 (\Omega_2) $.
  It is clear that $ w \in L^2 (\Omega) $; we only need to prove that $ \div \mu = w $.
  
  Let $ \varphi $ be a $ C^\infty $ function whose support is compact in $ \Omega $.
  By applying Green's formula separately in $ \Omega_1 $ and in $ \Omega_2 $ like in the proof of
  Lemma \ref{lemma:traces-3} and by noting that the integrals on $ \Gamma $ cancel out, we conclude that
  \begin{equation*}
    \int_\Omega w\, \varphi \,dx = -\int_{\Omega} \mu_i\, \partial_i \varphi\, dx
  \end{equation*}
  which means that $ \div \mu = w $.
\end{proof}

\begin{remark}\label{rem:traces-1}
  Recall that $ n^{(2)} = -n^{(1)} $.
  We may prefer to use only one normal vector; let $n$ denote either of $ n^{(1)} $ or $ n^{(2)} $.
  Then, the conclusion of Lemma \ref{lemma:traces-3} becomes\,: the two traces on $ \Gamma $ of the normal
  component of $ \mu $ are equal.
  The hypothesis of Lemma \ref{lemma:traces-4} becomes\,: $ tr^{(1)}_\Gamma (\mu^{(1)} n) =
  tr^{(2)}_\Gamma (\mu^{(2)} n) $.
\end{remark}
 

Lemmas \ref{lemma:traces-1}, \ref{lemma:traces-2}, \ref{lemma:traces-3} and \ref{lemma:traces-4}
can be found in \cite[Theorems 3.2.1 and 3.2.2]{ValliQu2008} for the context of finite elements.

\begin{figure}[ht] \centering
  \psfrag{A}{$A$}
  \psfrag{B}{$B$}
  \psfrag{C}{$C$}
  \psfrag{D}{$D$}
  \psfrag{Vin}{$ V_{in} $}
  \psfrag{Vout}{$ V_{out} $}
  \psfrag{Win}{$ W_{in} $}
  \psfrag{Wout}{$ W_{out} $}
  \includegraphics[width=55mm]{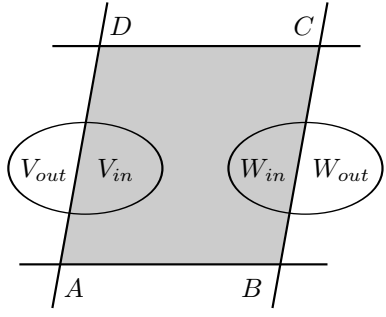}
  \caption{Periodicity cell}
  \label{fig cell bolas}
\end{figure}

Let us now go back to the periodic context.
Let $ u \in H^1_\# $; we shall focus on the behaviour of $u$ in a small open set $V$ crossing the side
$ AD $ of $ Y = ABCD $.
Let $ W = V + g_1 $ be the translation of $V$ by the periodicity vector $ g_1 = \vec{AB} = \vec{DC} $
and define $ V_{in} = V \cap Y $, $ W_{in} = W \cap Y $, $ V_{out} = V \setminus Y $,
$ W_{out} = W \setminus Y $ as in Figure \ref{fig cell bolas}.
Define also $ \Gamma_V = V \cap AD $ and $ \Gamma_W = W \cap BC $.

Since $ u \in H^1_\# $, $u$ has several restrictions\,: one in $ H^1(Y) $, nother one in
$ H^1(Y+g_1) $, another one on $ H^1(Y-g_1) $ and of course many other restrictions can be considered.
The restriction of $u$ on $Y$ has a trace on $ \Gamma_V $, the restriction of $u$ on $ Y-g_1 $ has
also its trace on $ \Gamma_V $.
Because $ u \in H^1(V) $, these two traces coincide (see Lemma \ref{lemma:traces-1}).
The same holds on $ \Gamma_W $.
Due to the periodic character of $u$, the trace on $ \Gamma_V $ coming from inside $Y$ is equal to
the trace on $ \Gamma_W $ coming from outside.
Similarly, the trace on $ \Gamma_W $ coming from inside $Y$ is equal to
the trace on $ \Gamma_V $ coming from outside.
Combining all the above (that is, combining the periodic character of $u$ with the fact that $u$ is
in $ H^1 $ of any ``small'' open set) we finally conclude the following.

\begin{lemma}\label{lemma:traces-5}
If $ u\in H^1_\# $, the traces of $u$ on opposite faces of $Y$ are equal.
\end{lemma}

The reverse process is also possible.
We can start with a function in $ H^1(Y) $ and extend it by periodicity;
the result will be a function in $ L^2_\# $ but not necessarily in $ H^1_\# $,
unless the traces of $u$ are compatible\,:

\begin{lemma}\label{lemma:traces-6}
  Let $ u\in H^1(Y) $ having equal traces on opposite faces of $Y$.
  Then $u$ can be extended by periodicity to a function in $ H^1_\# $.
\end{lemma}

Similar results hold for vector fields $ \mu $ in $ H_\# (\div) $.
Because the restriction of $ \mu $ to $Y$ belongs to $ H (\div,Y) $, its normal component has a
trace on $ \Gamma_V $.
Because the restriction of $ \mu $ to $ Y-g_1 $ belongs to $ H (\div,Y-g1) $, its normal component
has a trace on $ \Gamma_V $.
Because $ \mu\in H(\div,V) $, these two traces coincide (see Lemma \ref{lemma:traces-3} and
Remark \ref{rem:traces-1}).
Similarly, the two traces of the normal component of $ \mu $ on $ \Gamma_W $, one coming from
inside Y, the other one coming from outside, coincide.
On the other hand, the periodic character of $ \mu $ implies that the trace of its normal component
on $ \Gamma_V $ coming from inside $Y$ is equal to the trace on $ \Gamma_W $ coming from outside.
Similarly, the trace of the normal component of $ \mu $
on $ \Gamma_W $ coming from inside $Y$ is equal to the trace on $ \Gamma_V $ coming from outside.
Combining all the above (that is, combining the periodic character of $ \mu $ with the fact that
$ \mu $ is in $ H(\div,V) $ for any ``small'' open set V) we finally conclude the following.

\begin{lemma}\label{lemma:traces-7}
  If $ \mu\in H_\#(\div) $, the traces of the normal component of $ \mu $ on opposite faces of $Y$
  (in the spirit of Remark \ref{rem:traces-1}) are equal.
\end{lemma}

The reverse process is also possible.
We can start with a vector field $ \mu $ in $ H(\div,Y) $ and extend it by periodicity;
the result will be a field in $ L^2_\# $ but not necessarily in $ H_\#(\div) $,
unless the traces of the normal component of $ \mu $ are compatible\,:

\begin{lemma}\label{lemma:traces-8}
  Let $ \mu\in H(\div,Y) $ whose normal component has equal traces on opposite faces of $Y$,
  in the spirit of Remark \ref{rem:traces-1}.
  Then $ \mu $ can be extended by periodicity to a vector field in $ H_\#(\div) $.
\end{lemma}

Thus, in problem (\ref{eq:cell-pb}), the fact that $ \ffi $ is locally $ H^1 $ and is periodic
implies that $ \ffi $ has equal traces on opposite faces of $Y$;
the fact that $ \tens C\, \tens e(u_A) $ has zero divergence and is periodic
implies that the normal component of $ \tens C\, \tens e(u_A) $
has equal traces on opposite faces of $Y$, in the sense of Remark \ref{rem:traces-1}.
See also \cite[Section 2]{Luciano-Sacco}.

\section{By the way of compensated compactness}
\label{sec:compens-compact}

Theorem \ref{thrm:prod-avrg} below is a direct consequence of Green's formula.
But it also has far-reaching connections, as pointed out in Remark \ref{rem:hom-tens-energ-prod} below.
The equivalence between the variational formulations in
strain, in stress and in displacement in \cite{BT2022-B} is based on the Theorem \ref{thrm:prod-avrg}
which may be seen as a compensated compactness result, see Remark \ref{rem:comp-comp} below.

\begin{theorem}   \label{thrm:prod-avrg}
Given an arbitrary $ v \in LP $ and an arbitrary $ \tens\sigma \in \TT $, the following equality holds 
\begin{equation}  \label{eq:prod-avrg-sym}
\media_Y \langle \nabla\!{}_s v, \tens\sigma \rangle\, dy =
\bigl\langle \media_Y \nabla\!{}_s v\, dy, \media_Y \tens\sigma\, dy \bigr\rangle .
\end{equation}
\end{theorem}

\begin{proof}
  Define $ \displaystyle A = \media_Y \nabla\!{}_s v(y)\, dy $;
  then $ \ffi(y) = v(y) - Ay $ belongs to $ LP(0) = H^1_\#  $.
  Thus,
$$
\dint_Y \langle \nabla\!{}_s v , \tens\sigma \rangle\, dy =
\dint_Y \langle \nabla\!{}_s (A y + \ffi(y)) , \tens\sigma \rangle\, dy =
\langle A , \dint_Y \sigma\, dy \rangle  +
\dint_Y \langle \nabla\!{}_s \ffi , \tens\sigma \rangle\, dy \,.
$$
Green's formula (Theorem \ref{thrm:Green}) states that the second product is zero and this
concludes the proof.
\end{proof}
	 
\begin{remark} \label{rem:prod-avrg}
  For a given macroscopic stress $S$, the space $ \TT(S) $ is defined as the set of
  stress matrix fields in $ \TT $ whose mean is $S$.
  Similarly, for given macroscopic strain $A$, the space $ LP(A) $ is equal to the set of
  strain matrix fields in $ LP $ whose mean is $A$.
See \cite[Lemma 2]{BT2010}; the only difference is that in the present paper
the symmetry of the matrix $A$ is imposed in the very definition of the space $ LP $.

Thus, Theorem \ref{thrm:prod-avrg} above can be stated as follows.
  If $ v\in LP(A) $ and $ \tens\sigma \in \TT(S) $, 
\begin{equation} \label{eq:AS}
  \media_Y \langle \nabla\!{}_s v, \tens \sigma \rangle\, dy = \langle A, S \rangle .
\end{equation}
  This is the same as Proposition 1 stated by P.\ Suquet in \cite{Suquet1987}.
\end{remark}

\begin{remark} \label{rem:hom-tens-energ-prod}
  The above equality (\ref{eq:AS}) implies that the homogenized tensor $ \tens C^H $ can be defined
  as in formula (\ref{CHAB})
$$
\langle \tens C^H  A ,  B \rangle =
\bigl\langle \media_Y \tens C\, \tens e  (w_{ A}),  B \bigr\rangle =
\media_Y \langle \tens C\, \tens e  (w_{ A}),  B \rangle =
\media_Y \langle \tens C\, \tens e  (w_{ A}), \tens e  (w_{ B})\rangle ,
$$
where $ A$ and $ B$ are any two symmetric matrices (strains).
\end{remark}

The concept of compensated compactness, introduced by L.\ Tartar in \cite{Tartar1979},
describes the weak convergence of a product of two weakly convergent sequences
which must satisfy certain differential constraints.
It is also known as the div-curl lemma because in a particular case one of the sequences
has zero divergence while the other one has zero curl. We recall below the div-curl lemma as stated in \cite[Theorem 1]{Tartar1979}.

\begin{theorem}[Tartar's div-curl lemma] \label{div-curl-lemma}
  Let $ (a_n) $ and $ (b_n) $ be two sequences of vector fields in $ L^2(\Omega) $,
  weakly convergent towards $a$ and $b$, respectively.
  If the sequences $ div\, a_n $ and $ curl\, b_n $ belong to $ L^2(\Omega) $ and are
  bounded in that space, then the sequence of scalar products $ a_n \cdot b_n $ converges,
  in the sense of distributions, towards the scalar product $ a \cdot b $.
\end{theorem}

\begin{remark} \label{rem:comp-comp}
  There is a striking similarity between the hypotheses in Theorems \ref{div-curl-lemma} and
  \ref{thrm:prod-avrg}.
The similarity is somewhat hidden by the symmetric gradient; however, since $ \tens\sigma $
is a symmetric matrix, equation (\ref{eq:prod-avrg-sym}) can be re-written as
\begin{equation}  \label{eq:prod-avrg}
\media_Y \langle \nabla v, \tens\sigma \rangle\, dy =
\bigl\langle \media_Y \nabla v\, dy, \media_Y \tens\sigma\, dy \bigr\rangle
\end{equation}
In the above, we see that $ \nabla v $ has zero curl while $ \tens\sigma $ has zero divergence.

Thus, Theorem \ref{thrm:prod-avrg} above can be seen as a div-curl lemma in the periodic context.
Below we provide a proof of Theorem \ref{thrm:prod-avrg} based on the div-curl lemma.
\end{remark}

\begin{proof}[Alternative proof of Theorem \ref{thrm:prod-avrg}]
  Consider $v$ and $ \tens\sigma $ as in the statement of Theorem \ref{thrm:prod-avrg};
  we shall prove (\ref{eq:prod-avrg}) rather than (\ref{eq:prod-avrg-sym}).
  Choose any domain $ \Omega \subset \RR^3 $.
  For a fixed $i$, define $ a_n(x) $ as the vector field whose $j$th component is
  $ \tens\sigma_{ij} (x/n) $; define $ b_n(x) $ as the vector field whose
  $j$th component is $ \partial_j (v_i) = v_{i,j} (x/n) $.
  Then $ div\, a_n = 0 $ and $ curl\, b_n = 0 $ which is one of the requirements of
  Theorem \ref{div-curl-lemma}.
  A second requirement is the weak convergence of $ a_n $ and $ b_n $ in $ L^2 $;
  by Lemma \ref{lemma-oscillating}, they converge towards the constant vectors $a$ and $b$
  whose components are $ a_j = \displaystyle \media_Y \sigma_{ij}\, dy $ and
  $ b_j = \displaystyle\media_Y v_{i,j}\, dy $.
  Thus, Theorem \ref{div-curl-lemma} holds.
  On the other hand, Lemma \ref{lemma-oscillating} can be applied to the sequence of scalar products
  $ a_n \cdot b_n $ which thus converges, weakly in $ L^1 $, towards the constant scalar
  $ \displaystyle \media_Y \sum_j v_{i,j}, \sigma_{ij}\, dy $ (recall that $i$ is fixed).
  Combining these, we conclude the equality (\ref{eq:prod-avrg}).
\end{proof}

We end this section with a well-known result regarding oscillating functions,
whose proof can be found in \cite[Lemma A.1]{BM1983}.

\begin{lemma} \label{lemma-oscillating}
  Let $f$ be a periodic function whose restriction to the periodicity cell $Y$ belongs to $ L^p(Y) $,
  with $ p \ge 1 $.
  Define the sequence of oscillating functions $ f_n(x) = f (x/n) $. 
  Then the sequence $ f_n $ converges towards the constant scalar
  $ \displaystyle\media_Y f(y)\, dy $ weakly in $ L^p(\Omega) $, for any bounded domain $ \Omega $.
\end{lemma}

\section{Extensions of Donati's Theorem}
\label{sec:LP-T-Donati}

In the present Section we state and prove a series of extensions of Donati’s Theorem in the
periodic framework.
Besides their intrinsic importance, these results are useful for describing the duality relations
between variational formulations in displacement, in strain and in stress of the cellular problem
arising in the theory of homogenization, see
\cite{BT2022-B} which can be viewed as a sequel of the present paper.

We begin with a result that will be needed in the proof of the first extension of Donati's Theorem
(Theorem \ref{thrm:Donati.1} below)
and that is an adaptation to the periodic case of Theorem 3.3 in \cite{Amrouche2006}.

\begin{theorem}
\label{thrm:isomorphism}
The operator $\nabla\!{}_s : H^1_{\# 0} \mapsto \LL^2_{\# s}$ is an isomorphism from the space 
$$
H^1_{\# 0} : = \{ v \in  H^1_\# (\RR^3,\RR^3) : \media_Y v\, dy =0  \}
$$
onto $\Im \nabla\!{}_s $. Consequently $\Im \nabla\!{}_s $ is closed in $\LL^2_{\# s}$.

Moreover, the dual operator of $\nabla\!{}_s: H^1_{\# 0} \mapsto \LL^2_{\# s}$ is 
$-\div : \LL^2_{\# s} \mapsto H^{-1}_\#$.
\end{theorem}

\begin{proof}
Note that $\Ker \nabla\!{}_s = \{ 0 \}$ since the rotation is canceled by the periodicity assumtion and the translation is canceled by
the zero mean condition.
In order to prove that $\nabla\!{}_s$ is an isomorphism, we shall prove that there exists a positive constant $C$ such that for all displacement 
$v \in H^1_{\# 0}$
\begin{equation}
\label{eq:inequality_1}
\|v\|_{H^1_\#} \le C \| \nabla\!{}_s v \|_{\LL^2_\#}.
\end{equation}
The proof will be done by contradiction. Assume that the above inequality is false. 
Then, for each natural number $n \in \NN$ there exists $v_n \in H^1_{\# 0}$ such that
$$
\|v_n\|_{H^1_\#} > n \| \nabla\!{}_s v_n \|_{\LL^2_\#}.
$$
Consider the normalized functions $w_n:= \dfrac{v_n}{\|v_n\|_{H^1_\#}}$ that consequently verify, when $n \to \infty$,
$$
w_n \in H^1_{\# 0} , \quad \|w_n\|_{H^1_\#}= 1 \text{ and } \| \nabla\!{}_s w_n \|_{\LL^2_\#} \to 0 .
$$
By Rellich Theorem it turns out that there exists a subsequence $w_{n_k}$ of the sequence $w_n$
and a function $w\in L^2_\# = L^2(Y)$ such that
$$
w_{n_k} \to w \hbox{ in } L^2_\# .
$$
For the subsequence of the symmetric gradients holds
$$
\nabla\!{}_s w_{n_k} \to 0 \hbox{ in } \LL^2_\# .
$$
Then $(w_{n_k})$ is a Cauchy sequence with respect to the norm $\| \cdot \|_{L^2_\#} + \| \nabla\!{}_s \cdot \|_{\LL^2_\#}$ and by Korn's
inequality (\ref{eq:Korn_ineq}) it turns out that $(w_{n_k})$ is also a Cauchy sequence with respect to the norm $\| \cdot \|_{H^1_\#}$. Consequently
$w$ belongs to $H^1_{\# 0}$, hence it has zero mean and therefore $w=0$, which contradicts  $\|w\|_{H^1_\#}= 1$.
This proves inequality (\ref{eq:inequality_1}). The inverse inequality is obvious and therefore $\nabla\!{}_s$ is an isomorphism and thus 
the space of strain fileds $\Im \nabla\!{}_s$ is closed.

\noindent The last assertion is obtained from the Green's formula (\ref{eq:Green}) and it turns out that the dual operator of $\nabla\!{}_s$ is $-\div$.

\end{proof}

The following result is an extension of Donati's Theorem and adapts to the periodic context
the Theorem 4.2 in \cite{Amrouche2006}.
Under periodicity conditions we prove that a strain matrix field $\tens e$ is orthogonal
to all stress matrix fields $ \tens s $ with zero divergence if and only if it can be represented 
as the symmetric gradient of a displacement field $v $ having zero mean.

\begin{theorem} \label{thrm:Donati.1}
Consider $\tens e \in \LL^2_{\# s}$. 
Then there exists a unique $ v \in H^1_{\# 0}$ such that $ \tens e = \nabla\!{}_s v $ in $ \LL^2_{\# s} $
if and only if 
\begin{equation}
\dint_Y \langle \tens e , \tens s\rangle\, dy =0, \hbox{ for all } \tens s \in \LL^2_{\# s} \hbox{ such that } \div\, \tens s =0 \hbox{ in }
H^{-1}_\# .
\label{eq:Donati.1}
\end{equation}
\end{theorem}
	
\begin{proof}
Acording to Theorem \ref{thrm:isomorphism}, the space $\Im \nabla\!{}_s$ is closed in $\LL^2_{\# s}$ and the dual operator of 
$\nabla\!{}_s: H^1_{\# 0} \mapsto \LL^2_{\# s}$ is 
$-\div : \LL^2_{\# s} \mapsto H^1_{\# 0}$. Applying Banach's closed range theorem it turns out that $\Im \nabla\!{}_s = (\Ker\, \div )^\perp $ which
is equivalent to the claim in the theorem, that is, an element $\tens e \in \LL^2_{\# s}$ may be written as $\tens e = \nabla\!{}_s v$, for a function
$v \in H^1_{\# 0}$, iff $\tens e$ is orthogonal on all matrix functions $\tens s \in \LL^2_{\# s}$ with zero divergence $ \div\,\tens s =0 \hbox{ in }
H^{-1}_\#$. 
Since $\Ker \nabla\!{}_s = \{ 0 \}$, the function $v \in H^1_{\# 0}$ is unique.  

\end{proof}

\begin{remark} \label{rmk:Donati}
In the above Theorem \ref{thrm:Donati.1}, the average of the strain tensor $\tens e$ is zero\,:
$ \displaystyle\media_Y \tens e\, dy = 0 $.
This happens for two distinct reasons.
On one hand, if we take $ \tens s $ constant, it turns out from (\ref{eq:Donati.1}) that
$ \dint_Y \tens e\, dy = 0 $ and therefore the average is zero.
On the other hand, since $v$ belongs to $ H^1_{\# 0} $, $v$ is periodic
hence $\displaystyle\media_Y \nabla\!{}_s v\, dy = 0$.
\end{remark}

\begin{remark}
Theorem \ref{thrm:Donati.1} may be stated with $ v \in H^1_\# $. In this case $v$ is unique up to an additive constant vector.
If we add the zero-average hypothesis, $ v \in H^1_{\#0} $, then $v$ is unique.
\end{remark}

Recall the definition (\ref{eq:def-TT}) of $ \TT $;
recall that $ \TT(0) $ denotes the set of stresses belonging to $ \TT$ having zero average. 
Then the following Theorem is also an extension of Donati's Theorem,
where the test functions belong to $ \TT(0) $:
A strain matrix field $\tens e$ is orthogonal to all stress matrix fields $ \tens \mu $
with zero divergence and zero average if and only if it can be represented 
as the symmetric gradient of a displacement field $ v $ which is the sum between
a linear function and a periodic one.

\begin{theorem} \label{thrm:Donati.2}
Let $ \tens e \in \LL^2_{\# s} $. Then 
\begin{equation}
\dint_Y \langle \tens e, \tens\mu \rangle =0 \hbox{ for all } \tens\mu \in \TT(0)
\label{eq:Donati.2}
\end{equation}
if and only if there exists $ v \in LP $ such that $ \tens e= \nabla\!{}_s v $.
In this case $v$ is unique up to an additive constant vector.
\end{theorem}

\begin{proof}
Define $ A = \displaystyle\media_Y \tens e\, dy $ and consider
$ \tilde {\tens e} = \tens e - A $, thus $ \displaystyle\media_Y \tilde{\tens e}\, dy = 0 $.
The right hand member in (\ref{eq:Donati.2}) may be written as 
$$
\dint_Y \langle\tens e, \tens \mu \rangle = \dint_Y \langle \tilde{\tens e}, \tens \mu \rangle + 
\dint_Y \langle\displaystyle\media_Y \tens e, \tens \mu \rangle = 
\dint_Y \langle \tilde{\tens e}, \tens \mu \rangle + 
\langle \displaystyle\media_Y \tens e , \displaystyle\media_Y \mu  \rangle =
\dint_Y \langle\tilde{\tens e}, \tens \mu \rangle .
$$
Then the condition (\ref{eq:Donati.2}) becomes\,:
$$
\dint_Y \langle\tilde{\tens e}, \tens \mu \rangle =0 \hbox{ for all } \tens \mu \in \TT(0) .
$$
Theorem \ref{thrm:Donati.1} may be applied with $\tilde{\tens e}$ and consequently
the above condition is verified if and only if
there exists $\tilde v \in H^1_{\# 0}$ such that $ \tilde{\tens e} = \nabla\!{}_s \tilde v $.
Then $ \tens e = A + \nabla\!{}_s \tilde v = \nabla\!{}_s (A y +\tilde v) $.
Therefore we take $ v = Ay +\tilde v $ and thus $ v $ belongs to $ LP $.
\end{proof}

While in Theorem \ref{thrm:Donati.1} the mean of $ \tens e $ is zero (see Remark \ref{rmk:Donati}),
in Theorem \ref{thrm:Donati.2} the mean of $ \tens\mu $ is zero.
We can let both $ \tens e $ and $ \tens\mu $ to have non-zero average and state a result
(Theorem \ref{thrm:Donati.3} below) that
encompasses both Theorems \ref{thrm:Donati.1} and \ref{thrm:Donati.2}.
This result containts both Theorem \ref{thrm:prod-avrg} and its reciprocal.

\begin{theorem} \label{thrm:Donati.3}
Consider $\tens e \in \LL^2_{\# s}$. 
Then there exists $ v \in LP $ such that $ \tens e = \nabla\!{}_s v $ in $ \LL^2_{\# s} $
if and only if 
\begin{equation}
  \media_Y \langle \tens e , \tens s\rangle\, dy =
  \langle \media_Y \nabla\!{}_s v\, dy, \media_Y \tens s\, dy \rangle
  \hbox{ for all } \tens s \in \LL^2_{\# s} \hbox{ such that } \div\, \tens s = 0 .
\label{eq:Donati.3}
\end{equation}
\end{theorem}

\begin{proof}
  It suffices to introduce $ \tens\mu = \tens s - \displaystyle\media_Y \tens s $ and
  apply Theorem \ref{thrm:Donati.2}.
\end{proof}

In the following we state two results that give more information on the sets of stress matrix fields, displacement fields and strain matrix fields\,: 
the natural space $LP_0$ of displacement fields with zero average   is isomorphic to the orthogonal
complement of the stress matrix fields with zero divergence and zero average (Theorem \ref{thrm:LP-isomorphic-TT_0-perp})
and the space of periodic functions with zero average $H^1_{\# 0}$ is isomorphic to the orthogonal complement of the 
stress matrix fields with zero divergence, which is the natural space for strain matrix fields (Theorem \ref{thrm:H1-isomorphic-TT-perp}).
These results will allow one to
conclude the equivalence between the formulation in displacement and the formulation in strain as detailed in \cite{BT2022-B}.
Note that the following inclusions hold\,:
$$
H^1_{\# 0} \subset LP_0 \ \hbox{ and } \ \TT^\perp \subset \TT(0)^\perp \ \hbox{ while }
\ \TT(0) \subset \TT .
$$

\begin{theorem}[$H^1_{\# 0}$ is isomorphic to $\TT^\perp  $]
\label{thrm:H1-isomorphic-TT-perp}
The operator $\nabla\!{}_s : H^1_{\# 0} \mapsto \TT^\perp$ is an isomorphism.
\end{theorem}

\begin{proof}
In view of Theorem \ref{thrm:isomorphism}, $\Im \nabla\!{}_s$ is closed. The only point is to verify that $\Im \nabla\!{}_s = \TT^\perp$.
Consider a given displacement $v\in H^1_{\# 0}$ and an arbitrary stress field $\tens s\in \TT$. Denote $\tens e := \nabla\!{}_s v$. 
Green's formula (\ref{eq:Green}) 
implies that $\dint_Y \langle \tens e , \tens s \rangle d\ y =0$, therefore $\tens e$ belongs to $\TT^\perp$.  Conversely, taking an arbitrary 
$\tens e \in \TT^\perp$ and applying Theorem \ref{thrm:Donati.3} there exists $v\in LP$ such that  $\tens e = \nabla\!{}_s v$
and 
$
\displaystyle \media_Y \langle \tens e , \tens s\rangle\, dy =
  \langle \media_Y \nabla\!{}_s v\, dy, \media_Y \tens s\, dy \rangle =0
$
for all $\tens s \in \TT$. Hence $ \displaystyle\media_Y \nabla\!{}_s v =0$ so $v$ belongs to
$ H^1_{\# }$ and it may be chosen to belong to $H^1_{\# 0}$ and consequently it is unique.

\end{proof}

The following result states an isomorphism between the space $LP_0$ and the orthogonal complement of $\TT(0)$ by mean of the inverse operator 
of $\nabla\!{}_s$.

\begin{theorem}[$\TT(0)^\perp$ is isomorphic to $LP_0$]
\label{thrm:LP-isomorphic-TT_0-perp}
For each $\tens e \in \TT(0)^\perp $, denote by $\F(\tens e)$ the unique element in $LP_0$ that satisfies $ \nabla\!{}_s  \F(\tens e) = \tens e$
(according to Theorem \ref{thrm:Donati.2}). Then the mapping $\F : \TT(0)^\perp  \mapsto LP_0 $ defines an isomorphism between the Hilbert spaces
$\TT(0)^\perp$ and $LP_0$.
\end{theorem}

The proof is a consequence of Theorem \ref{thrm:Donati.2} and of the Korn inequality
in section \ref{sec:Korn-Green}.

\section*{Conclusions}
The present paper states and proves a series of extensions of Donati’s Theorem in the
periodic framework, in relation with the theory of homogenization.
This study allows for a mechanical interpretation of the relevant spaces\,:
a strain matrix field is orthogonal to the space of stress matrix fields with zero divergence
provided it can be represented as the symmetric gradient of a displacement field.
If we impose more conditions on the space of stress matrix fields the
corresponding strain and displacement fields adquire more properties.
These results are useful for obtaining equivalent variational formulations of the cellular problem
in stress, strain and displacement,
which are the object of another work by the same authors \cite{BT2022-B}.

A self-contained study of properties of traces in the periodic context was developed.
Besides being useful for proving Green's formula,
it allows for a better understanding of the behaviour of periodic functions and their traces.

\section*{Aknowledgements}

This work is supported by national funding from FCT - Foundation for Science and Technology (Portugal),
project UIDB/04561/2020.

The authors thank Lu\'isa Mascarenhas for suggesting bibliographical references about Lemma \ref{lemma-oscillating}.

\bibliography{biblio}

\end{document}